\documentclass[12pt]{amsart}
\usepackage{amscd}
\usepackage{graphicx}
\usepackage{hyperref}
\usepackage{times,rlepsf}
  \oddsidemargin 0in \theoremstyle{plain} \topmargin 0in
\newtheorem{thm}{Theorem}[section]

\newtheorem*{thm*}{Theorem}
\newtheorem*{prop*}{Proposition}
\newtheorem*{cor*}{Corollary}
\theoremstyle{definition}

\theoremstyle{remark}

\numberwithin{equation}{section}

\newcommand{\tam}{\mathbb Z}

\newcommand{\tuc}{\mathbb T^3}

\begin{document}

\title[Elliptic open books on torus bundles over the circle]
       {Elliptic open books on torus bundles over the circle}
\author{Tolga Etg\"{u}}
\address{Department of Mathematics , Ko\c{c} University, Sariyer, Istanbul 34450 TURKEY}
\email{tetgu@ku.edu.tr}

\thanks{The author is partially supported by a GEBIP grant of the Turkish Academy of Sciences
and a CAREER grant of the Scientific and Technological Research Council of Turkey .}%


\date{\today}

\begin{abstract}

As an application of the construction of open books on plumbed
3-manifolds, we construct elliptic open books on torus bundles over
the circle. In certain cases these open books are compatible with
Stein fillable contact structures and have minimal genus.

\end{abstract}

\maketitle


\section{Introduction}

Even though it has been known how to construct a contact structure
starting from an open book decomposition of a closed 3--manifold
since the work of Thurston and Winkelnkemper \cite{tw}, only after
Giroux's fundamental work \cite{gi} on the correspondence between
contact structures and open books the latter started to play an
essential role in contact topology. This correspondence becomes
one-to-one when certain equivalence relations are taken into
account. As a consequence, given a contact structure there are
infinitely many open books that correspond to it. This leads to the
problem of finding the ``simplest'' open book compatible with a
given contact structure. One way to interpret simplicity is in terms
of the genus of a page. Etnyre \cite{e} proved that whenever the
contact structure is overtwisted there is a corresponding planar
open book, i.e. an open book with page genus equal to 0, and he also
gave nontrivial obstructions to the existence of planar open books
for fillable contact structures (see \cite{oss} for other
obstructions). In an earlier work Ozbagci and the author \cite{eo}
described how to construct open books given a plumbing description
of a 3--manifold. The open books constructed this way coincide with
the ones obtained by using the methods given by Gay \cite{g}. As an
application of this construction and using the plumbing description
of torus bundles over the circle \cite{km} we explain how to
construct elliptic open books, i.e. open books with page genus equal
to 1, on torus bundles over the circle. We demonstrate this on many
examples and obtain the following results.

\begin{thm*}
There exists elliptic open books compatible with Stein fillable
contact structures on circle bundles over the torus with Euler
number less than $5$.
\end{thm*}

For negative Euler numbers the above statement was known before, and
for $\tuc$ as a result of Etnyre's work \cite{e} we knew that there
was no planar open book corresponding to the unique Stein fillable
contact structure. The following consequence of the above theorem
was previously announced by Van Horn \cite{vh}.

\begin{cor*}
The minimal page genus among all the open books compatible with the
unique Stein fillable contact structure on $\tuc$ is $1$.
\end{cor*}

There are exactly seven Seifert fibred spaces (besides the circle bundles over the torus) among torus bundles
over the circle. The next statement is about these 3--manifolds.

\begin{thm*}
On the Seifert fibred spaces $M(\pm\frac23,\mp\frac13, \mp\frac13)$,
$M(\pm\frac12,\mp\frac14,\mp \frac14)$, $M(\pm\frac12,\mp\frac13,
\mp \frac16)$ and $M(\frac12,\frac12,-\frac12,-\frac12)$ there exists elliptic open books  compatible 
with
Stein fillable contact structures.
\end{thm*}

Note that on many of the 3-manifolds considered in the above statements, there is a unique 
Stein fillable contact structure (up to contactomorphism, and in some cases up to contact isotopy). 
For this observation one combines the recent result of Gay \cite{gay} stating that contact structures
with positive Giroux torsion cannot be even strongly symplectically fillable with Giroux's classification of tight
contact structures on torus bundles over the circle specifically Theorem 1.3 in  \cite{giroux}  (cf. \cite{h})
which specifies the cases where the Giroux torsion uniquely determines the tight contact structure. 

\noindent {\bf Acknowledgement.} I would like to thank Burak Ozbagci
for helpful conversations and comments on a draft of this paper.

\section{Contact structures and open books}

We will assume throughout this paper that a contact structure
$\xi=\ker \alpha$ is coorientable (i.e., $\alpha$ is a global
1--form) and positive (i.e., $\alpha \wedge d\alpha >0 $ ). In the
following we describe the compatibility of an open book
decomposition with a given contact structure on a 3--manifold.

Suppose that for an oriented link $L$ in a 3--manifold $Y$ the
complement $Y\setminus L$ fibers over the circle as $\pi \colon Y
\setminus L \to S^1$ such that $\pi^{-1}(\theta) = \Sigma_\theta $
is the interior of a compact surface bounding $L$, for all $\theta
\in S^1$. Then $(L, \pi)$ is called an \emph{open book
decomposition} (or just an \emph{open book}) of $Y$. For each
$\theta \in S^1$, the surface $\Sigma_\theta$ is called a
\emph{page}, while $L$ the \emph{binding} of the open book. The
monodromy of the fibration $\pi$ is defined as the diffeomorphism of
a fixed page which is given by the first return map of a flow that
is transverse to the pages and meridional near the binding. The
isotopy class of this diffeomorphism is independent of the chosen
flow and we will refer to that as the \emph{monodromy} of the open
book decomposition.

An open book $(L, \pi)$ on a 3--manifold $Y$ is said to be
\emph{isomorphic} to an open book $(L^\prime, \pi^\prime)$ on a
3--manifold $Y^\prime$, if there is a diffeomorphism $f: (Y,L)  \to
(Y^\prime, L^\prime)$ such that $\pi^\prime \circ f = \pi$ on $Y
\setminus L$. In other words, an isomorphism of open books takes
binding to binding and pages to pages.

\begin{thm} [Alexander \cite{al}] Every closed and oriented
3--manifold admits an open book decomposition.
\end{thm}

In fact, every closed oriented 3--manifold admits a planar open
book, which means that a page is a disk $D^2$ with holes \cite{rol}.
On the other hand, every closed oriented 3--manifold admits a
contact structure \cite{marti}. So it seems natural to strengthen
the contact condition $\alpha \wedge d \alpha >0$ in the presence of
an open book decomposition on $Y$ by requiring that $\alpha>0 $ on
the binding and $d\alpha>0 $ on the pages.

\begin{def}\label{compatible} An open book decomposition of a 3--manifold
$Y$ and a contact structure $\xi$ on $Y$ are called
\emph{compatible} if $\xi$ can be represented by a contact form
$\alpha$ such that the binding is a transverse link, $d \alpha$ is a
symplectic form on every page and the orientation of the transverse
binding induced by $\alpha$ agrees with the boundary orientation of
the pages.
\end{def}

\begin{thm} [Giroux \cite{gi}] \label{giroux} Every contact
3--manifold admits a compatible open book. Moreover two contact
structures compatible with the same open book are isotopic.
\end{thm}

Under the light of this correspondence a natural question that comes
to mind is the following:

\noindent {\bf Question.} Given a contact structure, what is the
minimal (page) genus of all the compatible open books?

As Etnyre proved in \cite{e}, if the contact structure is
overtwisted, there is always a compatible planar open book. He also
proved that for fillable contact structures, there is an obstruction
for having a compatible planar open book (also see \cite{oss} for
other obstructions).

\begin{thm}[Etnyre \cite{e}]\label{obstruction} If $X$ is a symplectic
filling of a contact 3--manifold $(M,\xi)$ carried by a planar open
book, then $b^+_2(X)=b_2^0(X)=0$.
\end{thm}

At this point we should also note that the monodromy of a compatible
open book may be a good indicator of the fillability and tightness
of the contact structure.

\begin{thm}[Giroux \cite{gi}]\label{stein} A contact structure is Stein fillable if and
only if it is compatible with an open book whose monodromy can be
written as a product of right-handed Dehn twists.
\end{thm}

\begin{thm}[Honda-Kazez-Mati\'{c} \cite{hkm}]\label{veering} A contact structure is tight if and
only if the monodromy of any compatible open book is right-veering.
\end{thm}

An important consequence of Theorem~\ref{veering} which applies to
some examples we construct in the later sections is that if the
monodromy of an open book can be written as a product of Dehn twists
along disjoint curves and at least one of these Dehn twists is
left-handed and parallel to a boundary component, then the monodromy
is not right-veering and hence the open book is compatible with an
overtwisted contact structure.

\section{Kirby diagrams and plumbing graphs of torus bundles over the circle}

\begin{figure}[ht]

  \begin{center}

\relabelbox \small {\epsfxsize=5.5in
  \centerline{\epsfbox{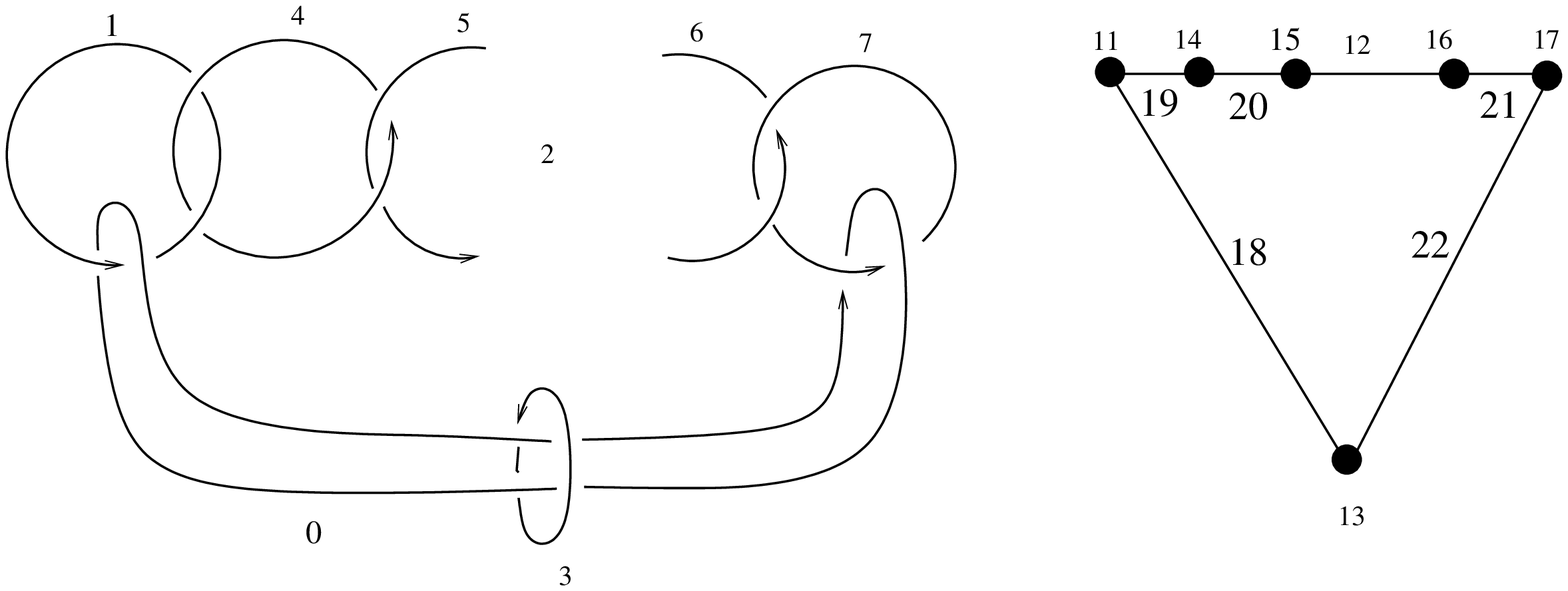}}}

  \relabel{0}{{$0$}}
  \relabel{1}{{$a_1$}}
  \relabel{2}{{$\cdots$}}
  \relabel{3}{{$0$}}
  \relabel{4}{{$a_2$}}
  \relabel{5}{{$a_3$}}
  \relabel{6}{{$a_{k-1}$}}
  \relabel{7}{{$a_k$}}
  \relabel{11}{{$a_1$}}
  \relabel{12}{{$\cdots$}}
  \relabel{13}{{$0$}}
  \relabel{14}{{$a_2$}}
  \relabel{15}{{$a_3$}}
  \relabel{16}{{$a_{k-1}$}}
  \relabel{17}{{$a_k$}}
  \relabel{18}{{$-$}}
  \relabel{19}{{$-$}}
  \relabel{20}{{$-$}}
  \relabel{21}{{$-$}}
  \relabel{22}{{$-$}}

   \endrelabelbox
   \caption{A framed link description of the torus bundle over the circle
    with monodromy $A=ST^{a_1}ST^{a_2}\cdots ST^{a_k}S$ and the corresponding plumbing graph}

\label{bundle}

    \end{center}

  \end{figure}

Let $Y$ be a 3-manifold which fibers over the circle with torus
fibers. It is uniquely determined by the monodromy of the fibration.
On the other hand, the automorphisms of the torus can be viewed as
(conjugacy classes of) matrices in $SL(2,\tam)$ by considering the
induced automorphisms on the first homology of the torus. In the
appendix to \cite{km}, Kirby and Melvin explain in detail how to
obtain a Kirby diagram on $Y$ given its monodromy $A$ as a matrix in
$SL(2,\tam)$. If $A=ST^{a_1}ST^{a_2}\cdots ST^{a_k}S$ for integers
$a_1,a_2,\dots , a_k$, where
$$S=\left(
      \begin{array}{cc}
        0 & -1 \\
        1 & 0 \\
      \end{array}
    \right)
\ \ \mbox{ and } \ \ T=\left(
                         \begin{array}{cc}
                           1 & 1 \\
                           0 & 1 \\
                         \end{array}
                       \right) $$
are matrices that generate $SL(2,\tam)$ with relations $S^4=I$ and
$S^2=(ST)^3$, then a Kirby diagram is given in Figure~\ref{bundle}
(left). As a consequence $Y$ has a plumbing description given by the
plumbing graph in Figure~\ref{bundle} (right). Note that the unknots in the Kirby diagram are oriented 
and according to these orientations the edges in the plumbing graph all have negative signs. If we had a tree as a plumbing graph the signs on the edges would be irrelevant as we could have changed them as we wish by changing the orientation of the base sphere of the circle bundle that corresponds to one of the two adjacent vertices. In our case we obviously don't have a tree and, even though we can change the signs on the edges in certain ways, we cannot, for example, make all the signs $+$ without changing the 3--manifold $Y$. 

\section{From plumbing diagrams to open books}


In \cite{eo}, Ozbagci and the author explain how to obtain open
books on 3--manifolds from a plumbing description of the manifold.
Even though the main focus of that paper is \emph{horizontal} open
books, a general procedure is described. Apparently, the open books
obtained this way coincide with the ones constructed by Gay
\cite{g}. Here we explain this construction on an example (this is
Example 10 in \cite{eo}).

\begin{figure}[ht]

  \begin{center}

     \includegraphics[scale=0.7]{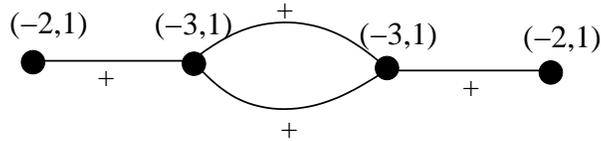}

   \caption{The integer weights at each vertex are the Euler number and the
genus of the base of the corresponding circle bundle, respectively.
Each edge is assumed to be positively signed.} \label{graph_ex1}

    \end{center}

  \end{figure}

Consider the 3--manifold $Y$ obtained by plumbing circle bundles
over tori according to the graph in Figure~\ref{graph_ex1}.

\begin{figure}[ht]

  \begin{center}

     \includegraphics{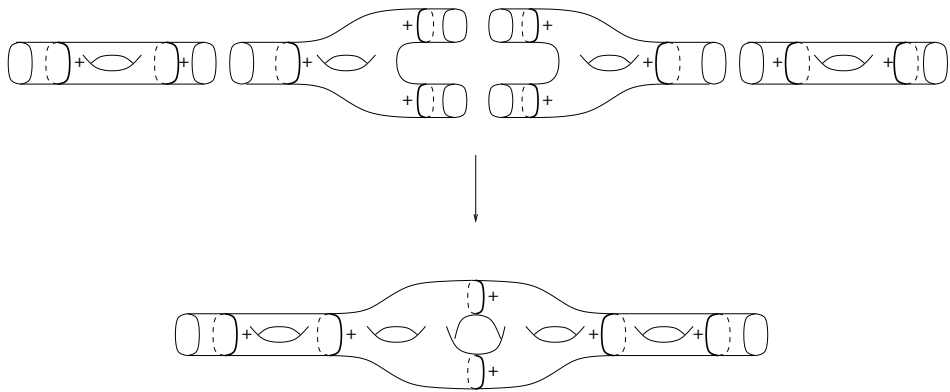}

   \caption{A page and the monodromy of the open book constructed from the plumbing above.} 
\label
{ex1_open_book}

    \end{center}

  \end{figure}

First of all, we construct an open book for each vertex, i.e. circle
bundle over a surface (these are shown in the upper part of Figure
~\ref{ex1_open_book}). In general for each vertex we obtain an open
book with a page of genus and punctures given by the base genus and
the absolute value of the Euler number, respectively, of the circle
bundle which corresponds to the vertex and the monodromy is the
product of Dehn twists around the punctures (right-handed if the
Euler number is negative and left-handed otherwise). Note that,
whenever necessary, we can add two boundary components to a page of
an open book without changing the underlying 3--manifold by adding
opposite-handed Dehn twists around these punctures. Going back to
our example, we then ``join'' the pages of these open books
according to the plumbing graph and keeping in mind that the edges
are all positive (see the lower part of Figure
~\ref{ex1_open_book}). In general, if there is a negative edge, to
join the pages of the open books that correspond to the adjacent
vertices we use boundary components that have left-handed Dehn
twists around them.

\section{Open books on torus bundles over the circle}

Combining the procedures described in the previous two sections, one
can obtain an elliptic open book on any 3--manifold which fibers
over the circle with torus fibers. Since we are interested in
realizing the minimal genus for a contact structure and since it is
already known that the minimal genus of any overtwisted contact
structures is zero \cite{e}, we should focus on tight contact
structures. On the other hand, even though we can detect the
tightness of a contact structure by considering open books
compatible with it \cite{hkm}, at least for the moment, we have to
consider \emph{all} the compatible open books. This is why we try to
construct open books compatible with Stein fillable contact
structures which are much easier to detect by considering the
monodromy of a \emph{single} compatible open book. On the other
hand, obtaining open books compatible with Stein fillable contact
structures requires some care. In what follows we focus on two
special classes among torus bundles over the circle: Seifert
fibred spaces and circle bundles over the torus.

\subsection{Open books on some Seifert fibred spaces}

\begin{figure}[ht]

  \begin{center}

     \includegraphics[scale=.5]{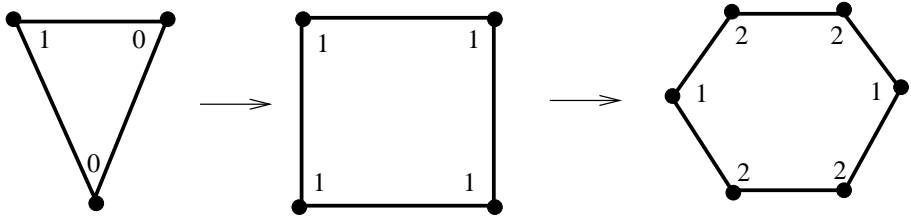}

   \caption{From the first graph to the second a single blow-up on the edge connecting 0-vertices and from the second to the third one blow-up on each side edge are performed. Every edge has a minus sign and every blow-up is with a $+1$-framed unknot.}

\label{blowups}

    \end{center}

  \end{figure}

As a first example, let us consider the small Seifert fibred
space $M(-\frac23,\frac13, \frac13)$. Using the fact that this
manifold fibers over the circle with monodromy
$$\left(
     \begin{array}{cc}
        0 & 1 \\
        -1 & -1 \\
     \end{array}
    \right) 
    =
    ST^1ST^0S
       $$
(see, for example, the  table on page 89 in \cite{h}), we provide an
open book decomposition of this manifold and the corresponding
plumbing graph in the first row of Table~\ref{listSF1}. Figure~\ref{blowups}
demonstrates how we obtain this plumbing graph from the initial plumbing graph via blow-ups.

\begin{table}
  \centering
  \caption{Elliptic open books on some Seifert fibred spaces.
The edges in the plumbing graphs have negative sign. Punctured tori
(black dots on the tori are the punctures) represent a page, the
monodromy is the product of Dehn twists along with the curves drawn
on this page and the handedness of these twists are given by the
sign next to each curve: $+$ indicates a right-handed (positive)
Dehn twist and $-$ indicates a left-handed (negative) Dehn
twist.}\label{listSF1}
\begin{tabular}{|c|c|c|c|}
  \hline
     3--manifold & monodromy & plumbing & open book
  \\ \hline
  \raisebox{1.5cm}[0pt]{$M(-\frac23,\frac13,\frac13) $} & \raisebox{1.5cm}[0pt]{$\left(
     \begin{array}{cc}
        0 & 1 \\
        -1 & -1 \\
     \end{array}
    \right)$}& \includegraphics[scale=.65]{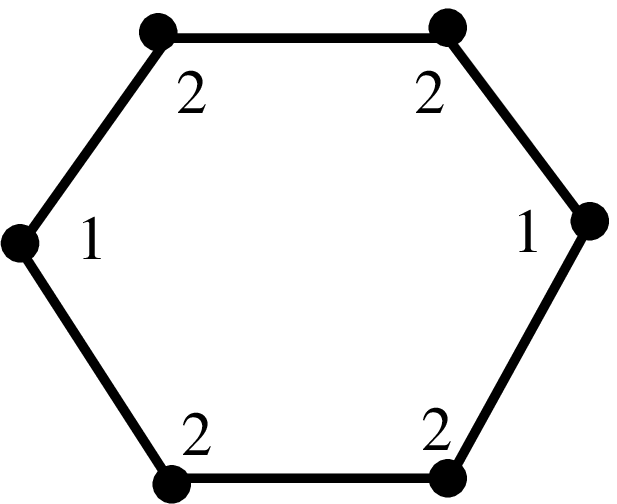} & \includegraphics[scale=.55]
{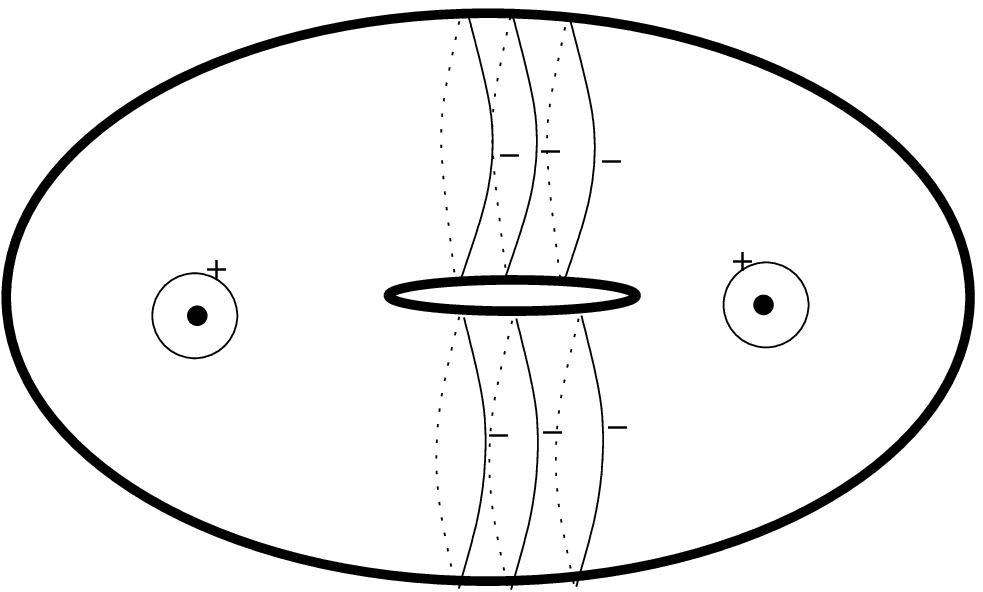}
  \\ \hline
   \raisebox{1.5cm}[0pt]{$M(-\frac12,\frac14,\frac14)$}& \raisebox{1.5cm}[0pt]{$\left(
     \begin{array}{cc}
        0 & 1 \\
        -1 & 0 \\
     \end{array}
    \right)$}&\includegraphics[scale=.65]{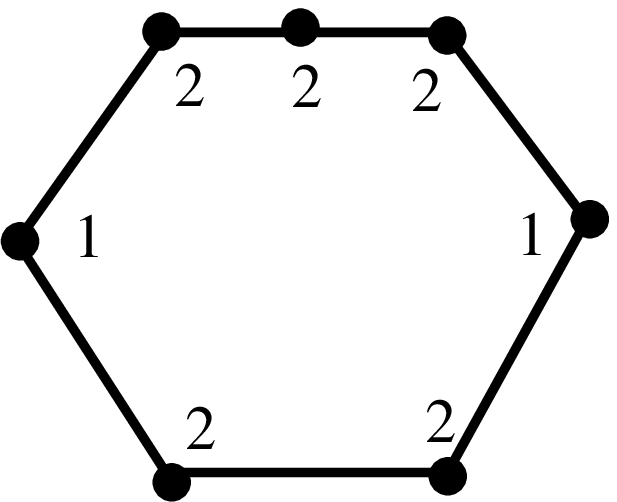} & \includegraphics[scale=.55]
{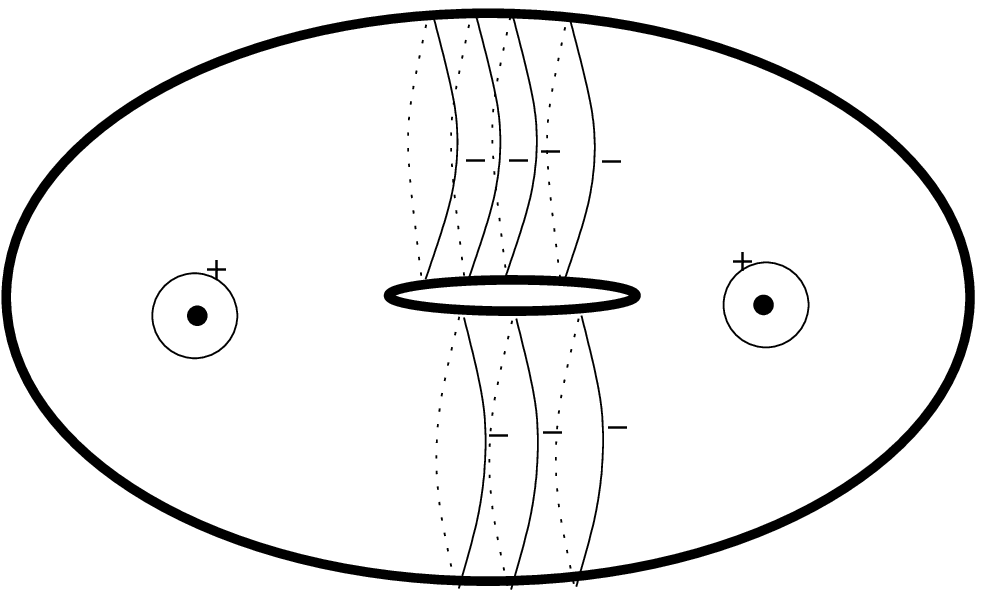}
   \\ \hline
  \raisebox{1.5cm}[0pt]{$M(-\frac12,\frac13,\frac16)$} &\raisebox{1.5cm}[0pt]{$\left(
     \begin{array}{cc}
        1 & 1 \\
        -1 & 0 \\
     \end{array}
    \right)$} &\includegraphics[scale=.65]{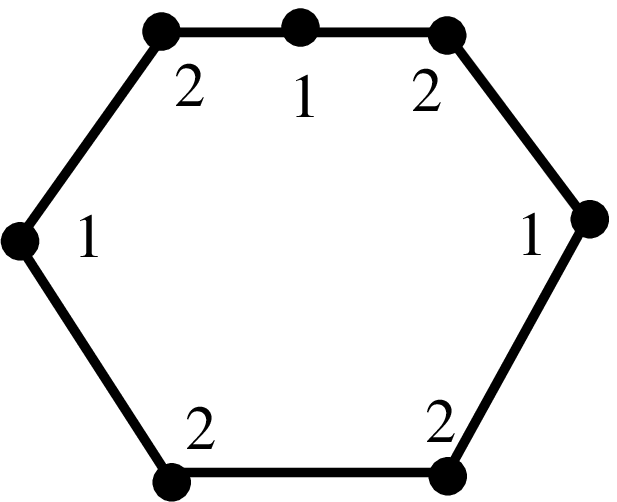} & \includegraphics[scale=.55]
{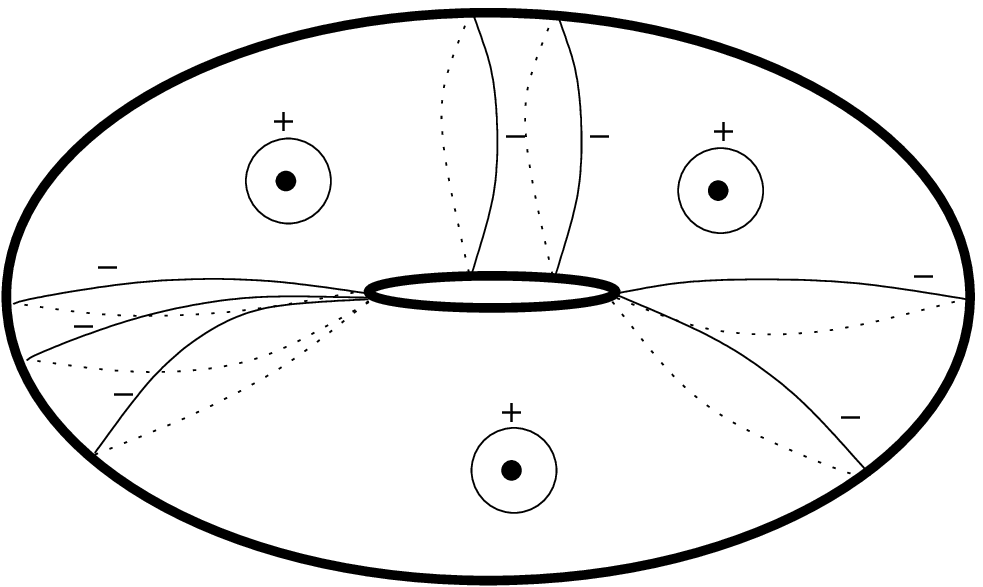}
  \\ \hline
  \raisebox{1.5cm}[0pt]{$M(-\frac12,-\frac12,\frac12,\frac12)$}& \raisebox{1.5cm}[0pt]{$\left(
     \begin{array}{cc}
        -1 & 0 \\
        0 & -1 \\
     \end{array}
    \right)$}& \includegraphics[scale=.65]{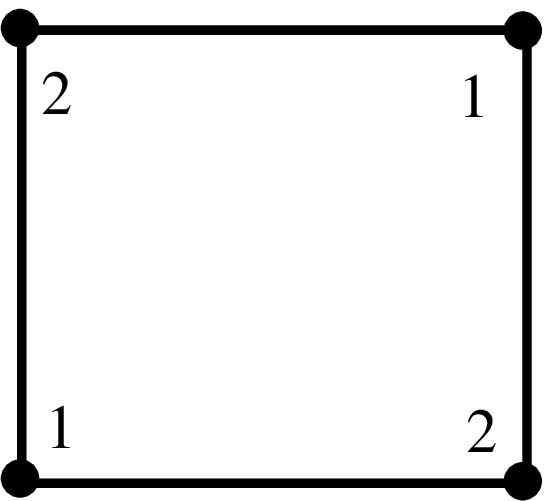} & \includegraphics
[scale=.55]{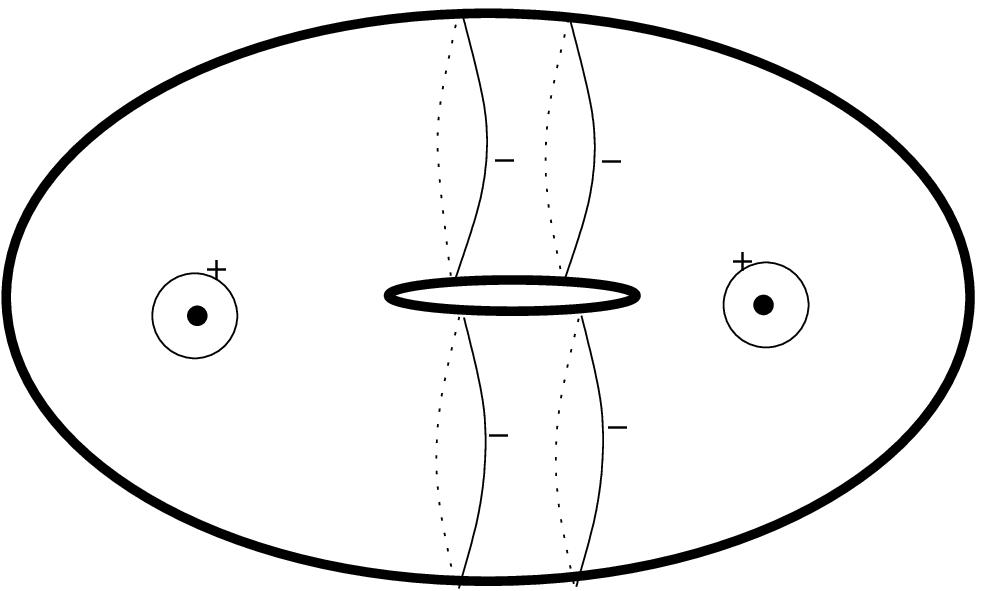}
  \\ \hline

\end{tabular}

\end{table}

\begin{figure}[ht]

  \begin{center}

\relabelbox \small {\epsfxsize=2.5in
  \centerline{\epsfbox{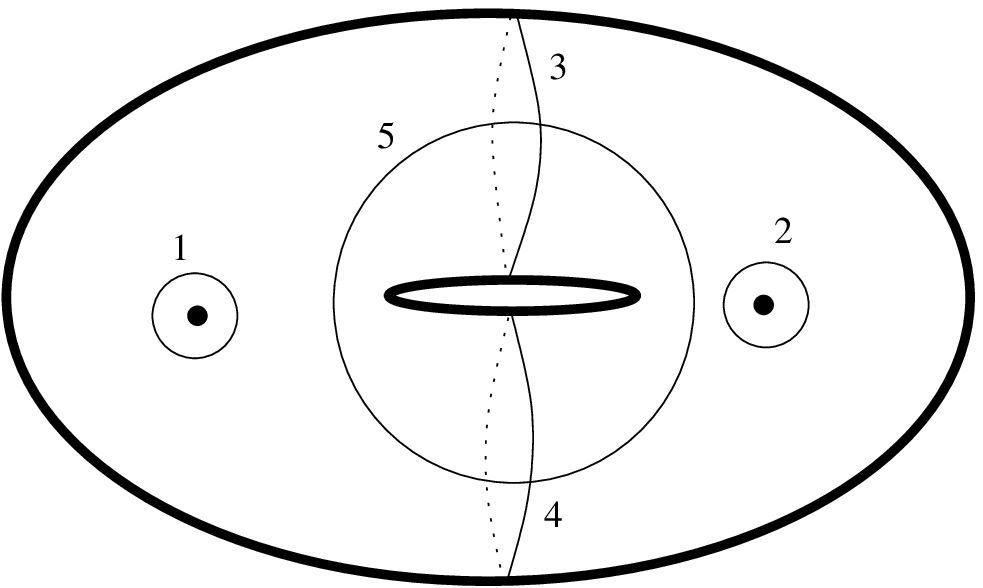}}}

  \relabel{1}{{$\delta_1$}}
  \relabel{2}{{$\delta_2$}}
  \relabel{3}{{$\alpha_1$}}
  \relabel{4}{{$\alpha_2$}}
  \relabel{5}{{$\beta$}}

   \endrelabelbox
   \caption{}

\label{twists}

    \end{center}

  \end{figure}

This open book is compatible with a Stein fillable contact structure
since its monodromy can be written as a product of right-handed Dehn
twists. To see this, let us denote the right-handed Dehn twists
around the punctures by $\delta_i$ and the right-handed Dehn twists
around the essential curves by $\alpha_i$ as in Figure~\ref{twists}.
So the monodromy is $\phi =
\delta_1\delta_2\alpha_1^{-3}\alpha^{-3}_2$ (note that since the
curves are pairwise disjoint, all these Dehn twists commute with
each other). We have the relation
$\delta_1\delta_2=(\alpha_1\alpha_2\beta)^4$, where $\beta$ is the
right-handed Dehn twists around the essential curve disjoint from
the $\delta_i$'s and intersecting each $\alpha_i$ geometrically once
(see \cite{ko} for this and other such relations on punctured tori).
It is crucial that we have no more than 4 parallel left-handed Dehn
twists on a twice-punctured torus in order to be able to use this
relation to write the monodromy in terms of right-handed Dehn
twists.

\begin{thm}
On the Seifert fibred space $M(-\frac23,\frac13, \frac13)$  there is
an elliptic open book compatible with a Stein fillable contact
structure.
\end{thm}

\begin{table}
  \centering
  \caption{Elliptic open books on some Seifert fibred spaces.}\label{listSF2}
\begin{tabular}{|c|c|c|c|}
  \hline
     3--manifold & monodromy & plumbing & open book
  \\ \hline

  \raisebox{1.5cm}[0pt]{$M(\frac23,-\frac13,-\frac13)$} &\raisebox{1.5cm}[0pt]{$\left(
     \begin{array}{cc}
        -1 & -1 \\
        1 & 0 \\
     \end{array}
    \right)$}& \includegraphics[scale=.65]{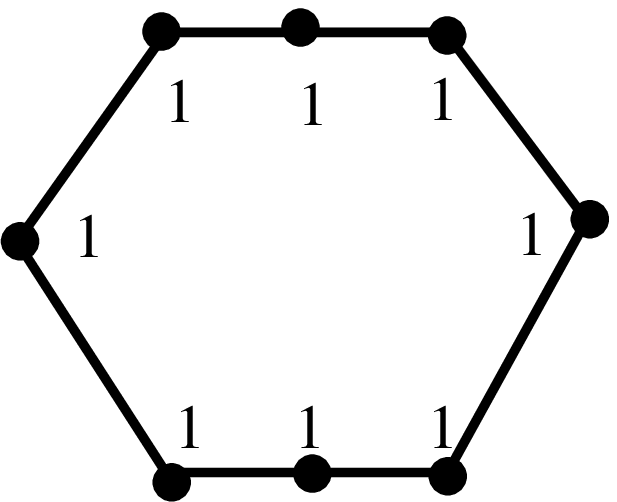} &\includegraphics[scale=.55]
{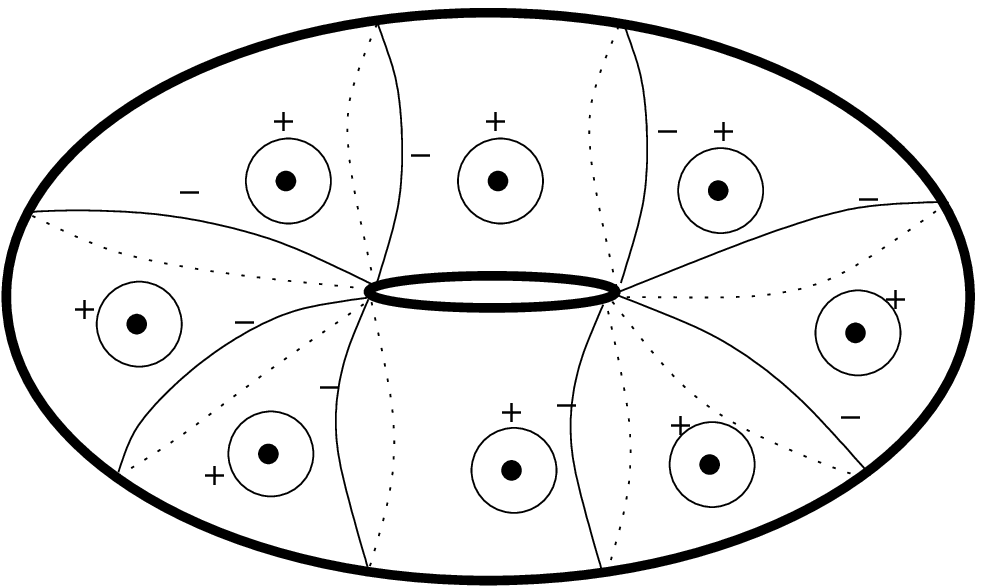}
  \\ \hline
   \raisebox{1.5cm}[0pt]{$M(\frac12,-\frac14,-\frac14)$} &\raisebox{1.5cm}[0pt]{$\left(
     \begin{array}{cc}
        0 & -1 \\
        1 & 0 \\
     \end{array}
    \right)$} &\includegraphics[scale=.65]{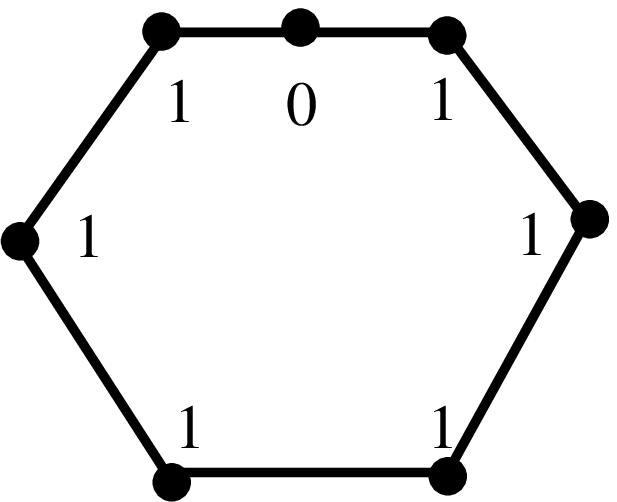} &\includegraphics[scale=.55]
{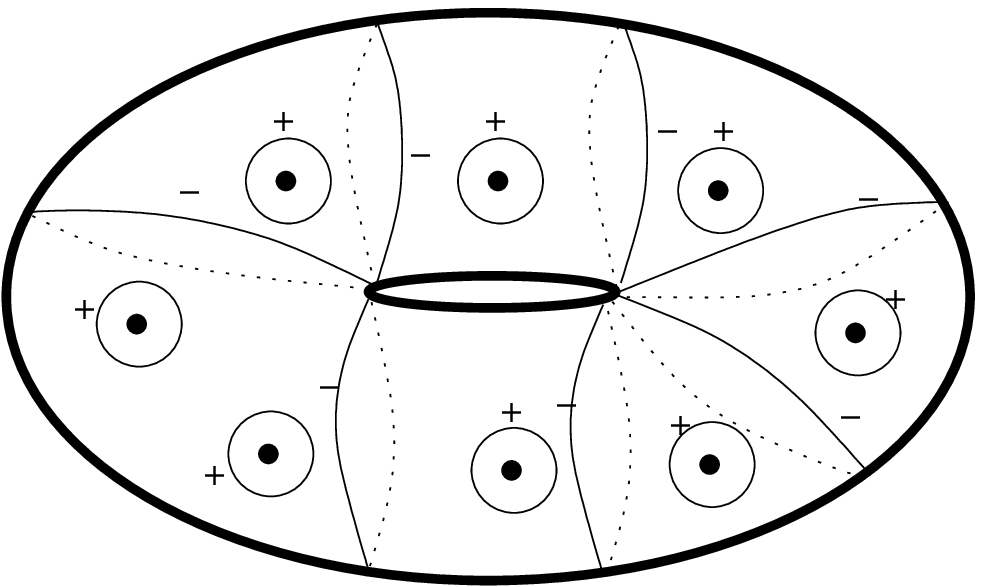}
  \\ \hline
   \raisebox{1.5cm}[0pt]{$M(\frac12,-\frac13,-\frac16)$}&\raisebox{1.5cm}[0pt]{$\left(
     \begin{array}{cc}
        0 & -1 \\
        1 &  1 \\
     \end{array}
    \right)$} &\includegraphics[scale=.65]{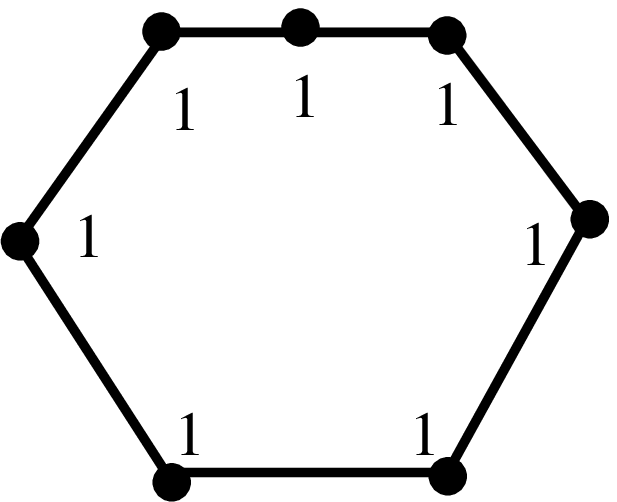} &\includegraphics[scale=.55]
{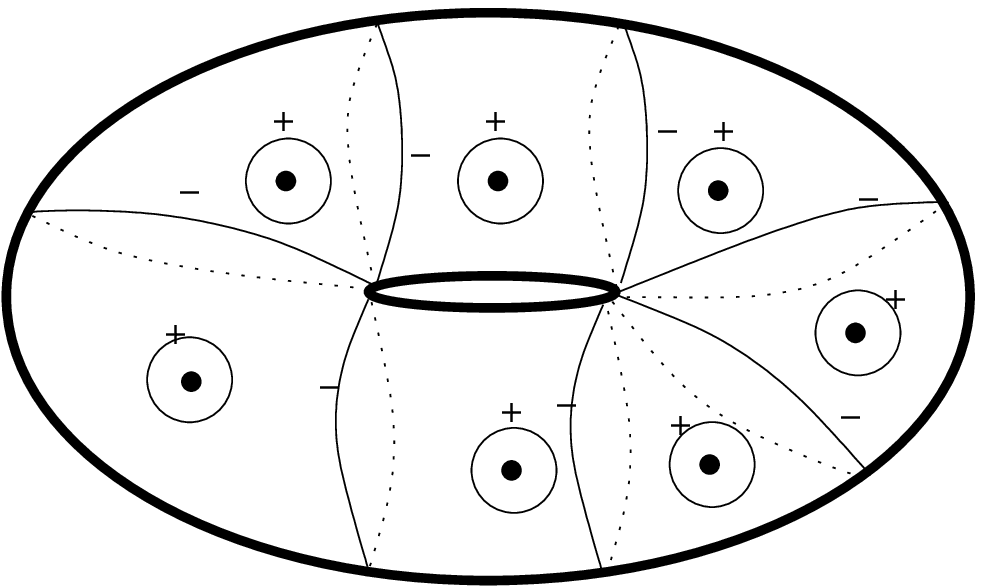}
  \\ \hline

\end{tabular}

\end{table}

In fact, using the same methods it is not difficult to do the same
for the other six Seifert fibred spaces which are torus bundles over
the circle. Examples of such open books are given in
Tables~\ref{listSF1} and \ref{listSF2}.

\begin{thm}
On the Seifert fibred spaces $M(\frac23,-\frac13, -\frac13)$,
$M(\pm\frac12,\mp\frac14, \mp\frac14)$, $M(\pm\frac12,\mp\frac13,
\mp \frac16)$, and $M(\frac12,\frac12,-\frac12,-\frac12)$ there are
elliptic open books compatible with Stein fillable contact
structures.
\end{thm}

\subsection{Open books on circle bundles over the torus}

In this subsection, we will give examples of open books on some
circle bundles over the torus. $Y_k$ denotes the circle bundle with
Euler number $k$. Since $Y_k$ is also a torus bundle over the
circle, we can use the previous techniques to obtain open books. As
in the previous subsection, we are especially interested in open
books compatible with Stein fillable contact structures.

\begin{figure}[ht]

  \begin{center}

     \includegraphics[scale=.45]{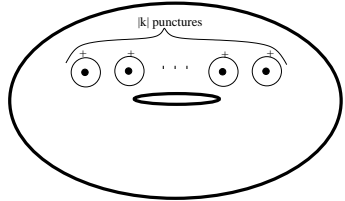}

   \caption{A page and the monodromy of the open book $\mathsf{ob}_k$ for $k<0$}

\label{ob_k}

    \end{center}

  \end{figure}

First of all, in case $k < 0$ (resp. $k>0$), as it was explained in
\cite{eo}, one can easily obtain an open book $\mathsf{ob}_k$ with
$|k|$-times punctured torus page and the product of one right-handed
(resp. left-handed) Dehn twist around each boundary component as
monodromy on $Y_k$. Note that when $k<0$, since the monodromy of
$\mathsf{ob}_k$ is the product of right-handed Dehn twists,
$\mathsf{ob}_k$ is compatible with a Stein fillable contact
structure on $Y_k$ (see Figure~\ref{ob_k}).

One obvious way to obtain an elliptic open book on $Y_0=\tuc$ is to
let a page be the twice punctured torus and the monodromy be the
product of opposite-handed Dehn twists around the punctures. But,
like $\mathsf{ob}_k$ for $k>0$, this open book is compatible with an
overtwisted contact structure since the monodromy is not right
veering \cite{hkm}.

\begin{figure}[ht]

  \begin{center}

     \includegraphics[scale=.75]{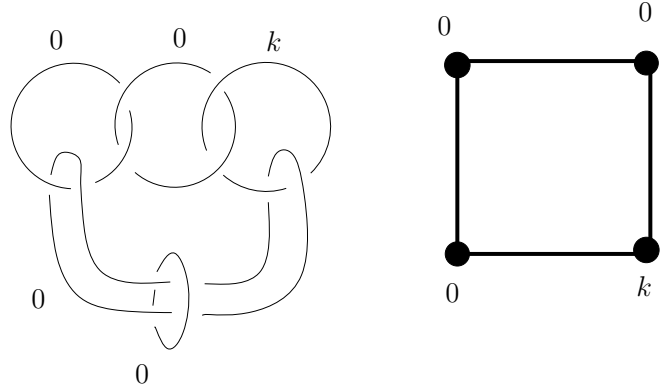}

   \caption{A framed link description of $Y_k$ and the corresponding plumbing graph}

\label{y_k}

    \end{center}

  \end{figure}

The total space $Y_k$ of the circle bundle of Euler number $k$ over
the torus is also the total space of the torus bundle over the
circle with monodromy $T^k = S^2T^kS^2 = ST^0ST^kST^0S$. Therefore
$Y_k$ can be given by the Kirby diagram and the corresponding
plumbing graph in Figure~\ref{y_k}.

\begin{figure}[ht]

  \begin{center}

     \includegraphics[scale=.55]{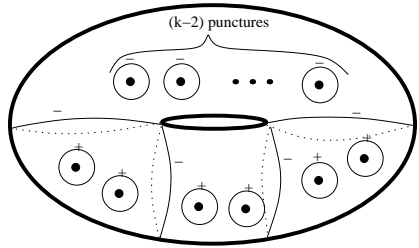}

   \caption{A page and the monodromy of the open book $\mathsf{ob}'_k$ for $k>2$}

\label{ob'_k}

    \end{center}

  \end{figure}

Using the recipe in \cite{eo} to get an open book on $Y_k$ from this
plumbing description we obtain an open book $\mathsf{ob}'_k$ with
$6+|k-2|$-times punctured torus page and a monodromy given by the
product of Dehn twists around the curves indicated in
Figure~\ref{ob'_k} for $k>2$.
Note that in case $k>2$ we get left-handed Dehn twists around $k-2$
boundary components (and the other Dehn twists are around disjoint
curves) which implies that the monodromy is not right-veering hence
the open book $\mathsf{ob}'_k$ is compatible with an overtwisted
contact structure \cite{hkm}. On the other hand, when $k\leq 2 $, we
have right-handed Dehn twists around boundary components. In
particular, when $k=0$ or $1$, by using certain relations \cite{ko}
in the mapping class group of torus with 8 and 7 punctures,
respectively , the monodromy can be written as a product of
right-handed Dehn twists implying that $\mathsf{ob}'_k$ is
compatible with a Stein fillable contact structure. Elliptic open
books compatible with the unique Stein fillable contact structure on
$\tuc$ are especially interesting as minimal genus open books
compatible with this contact structure since there is no planar open
book compatible with it: If there were a planar open book compatible
with the unique Stein fillable contact structure on $\tuc$, then by
Theorem 4.1 in \cite{e} any such filling would have vanishing
$b_2^+$ and $b_2^0$ while on the other hand $T^2 \times D^2$ has
$b_2^0\neq0$ and it is a Stein filling of $\tuc$.

\begin{thm}
The minimal page genus among all the open books compatible with the
unique Stein fillable contact structure on $\tuc$ is $1$.

\end{thm}

For $k=2$ case we weren't able to verify that the monodromy of
$\mathsf{ob}'_2$ can be written as a product of right-handed Dehn
twists, but obtained another elliptic open book on $Y_2$ which is
compatible with a Stein fillable contact structure. To see this,
take the plumbing diagram for $Y_2$ in Figure~\ref{y_k}. After
``blowing-up" this diagram four times we get the plumbing diagram and
the corresponding open book in the third row of Table~\ref{listY0}.
\begin{table}
  \centering
  \caption{Examples of elliptic open books on some circle bundles over the torus}\label{listY0}
\begin{tabular}{|c|c|c|}
  \hline
     3--manifold & plumbing & open book
  \\ \hline

  $Y_0=\tuc$& \includegraphics[scale=.65]{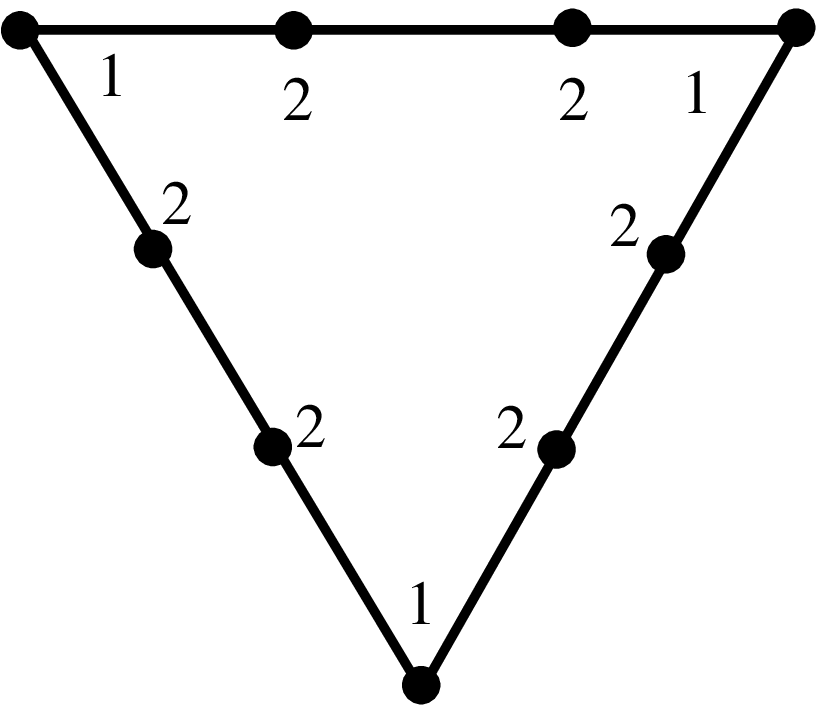} &\includegraphics[scale=.6]{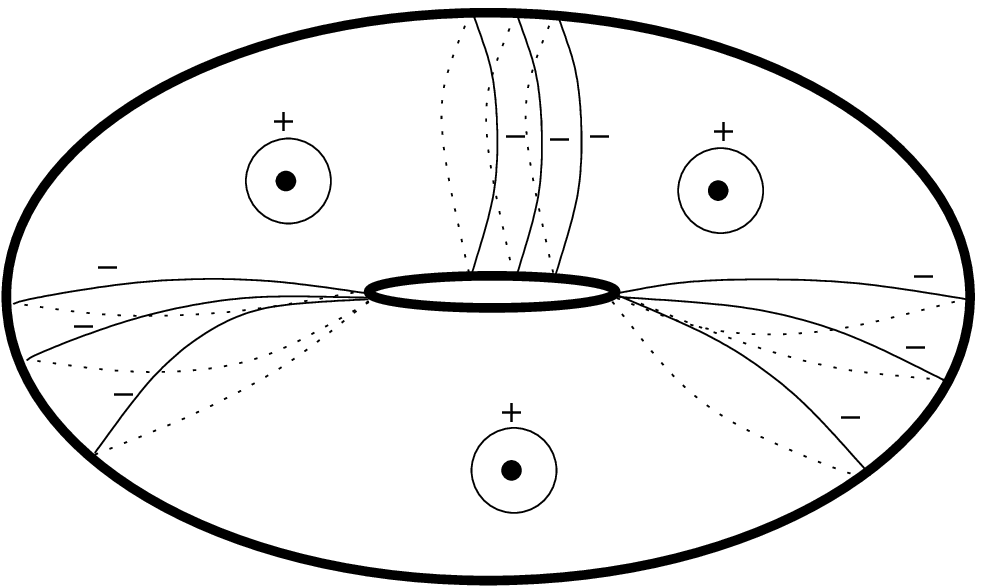}
  \\ \hline
  $Y_1$ & \includegraphics[scale=.65]{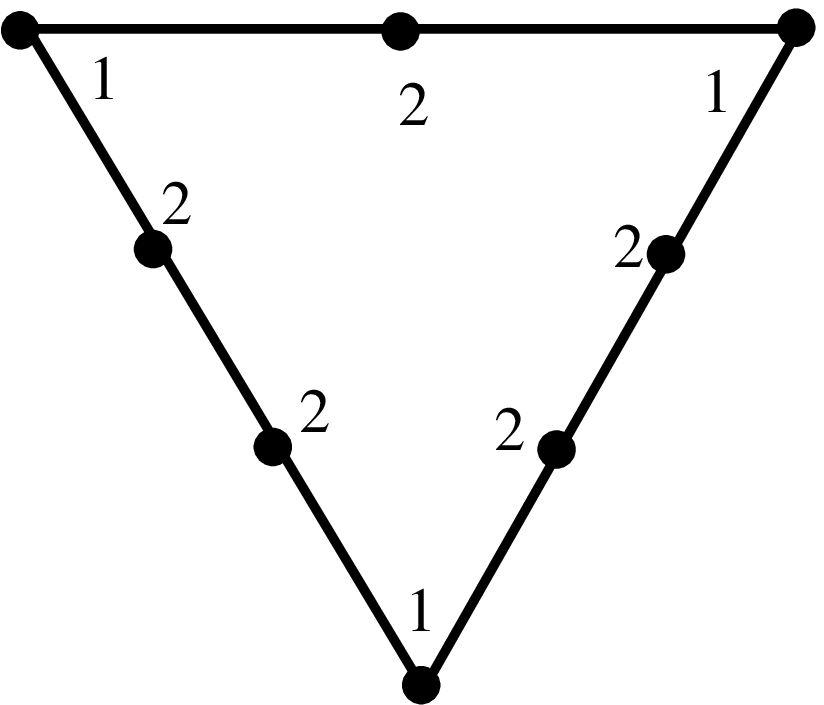} &\includegraphics[scale=.6]{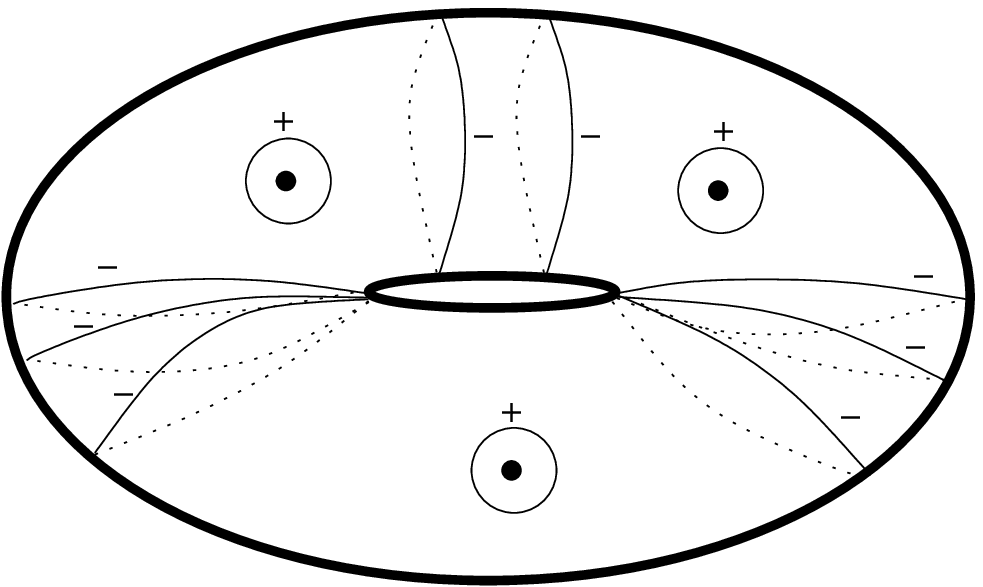}
   \\ \hline
  $Y_2$ & \includegraphics[scale=.6]{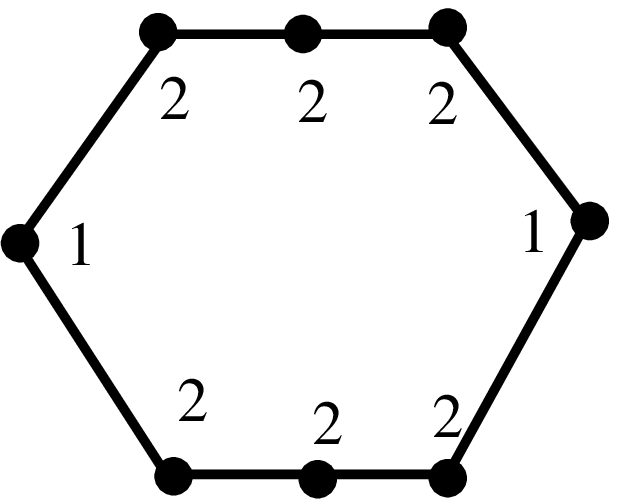} &\includegraphics[scale=.65]{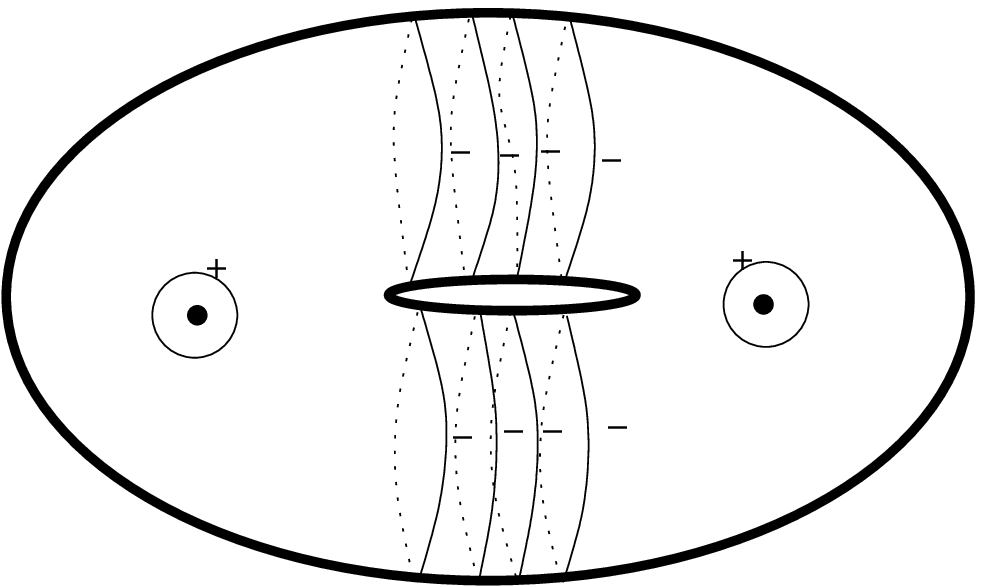}
  \\ \hline
  $Y_3$ & \includegraphics[scale=.65]{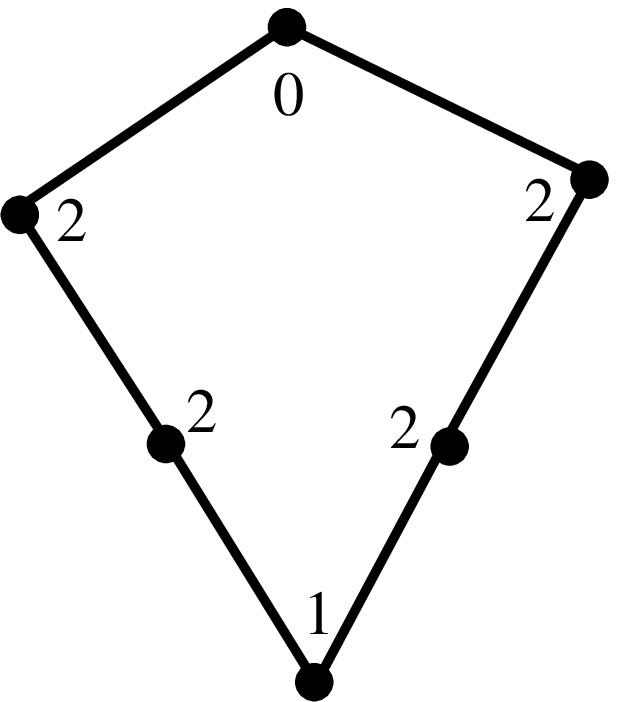} & \includegraphics[scale=.55]{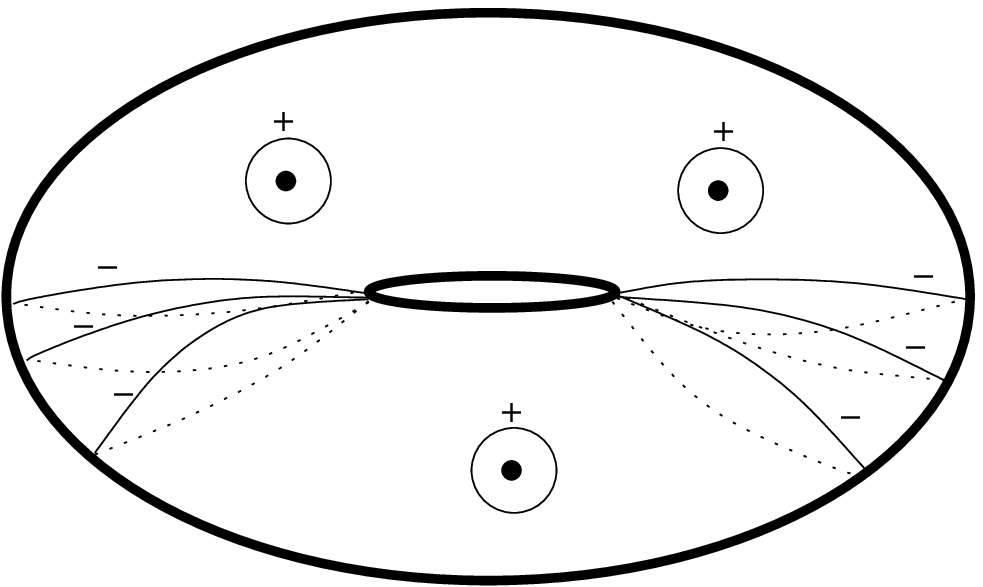}
  \\ \hline
  $Y_4$ & \includegraphics[scale=.65]{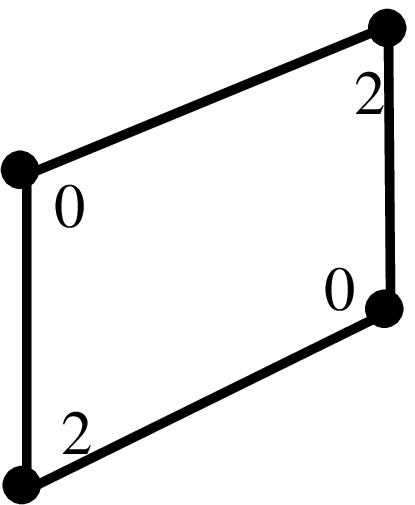} &\includegraphics[scale=.55]{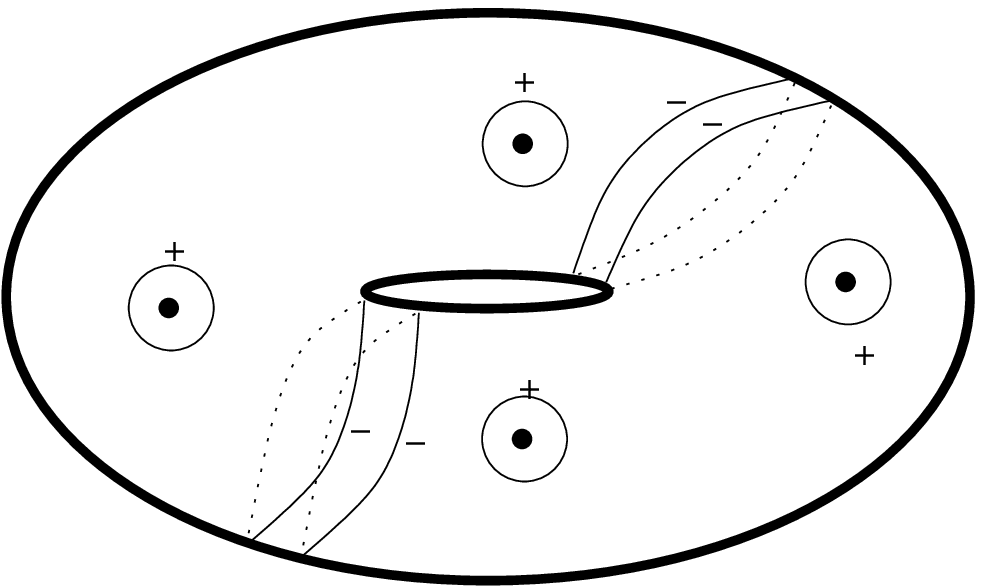}
  \\ \hline
\end{tabular}

\end{table}
Using a relation \cite{ko} in the mapping class group of the torus
with four punctures one can see that the monodromy of this new open
book can be written as a product of right-handed Dehn twists. In
fact, all of the open books given in Table~\ref{listY0} are obtained
this way, i.e. by blowing-up the plumbing diagram for $Y_k$ in
Figure~\ref{y_k} in a way that the monodromy of the resulting open
book can be written in terms of right-handed Dehn twists by using
the relations in \cite{ko}. We should note that the elliptic open
book given in the first row of Table~\ref{listY0} is compatible with the unique Stein fillable contact structure on
$Y_0=\tuc$ and it was first constructed
by Van Horn \cite{vh} using different techniques.

\begin{thm}
For $k \leq 4$ there exists elliptic open books on $Y_k$ compatible
with Stein fillable contact structures.
\end{thm}


\bibliographystyle{amsalpha}

\begin{thebibliography}{1}

\bibitem{al}
J. Alexander, \emph{A lemma on systems of knotted curves,} Proc.
Nat. Acad. Sci. USA \textbf{9} (1923), 93--95.

\bibitem{e}
J.~Etnyre, \emph{Planar open books and contact structures}, Int.
Math. Res. Not. (2004), 4255--4267.



\bibitem{eo}
T.~Etg\"{u} and B.~Ozbagci, \emph{Explicit horizontal open books on
some plumbings}, Internat. J. Math. \textbf{17} (2006), 1013--1031.

\bibitem{g}
D.~Gay, \emph{Open books and configurations of symplectic surfaces},
Algebr. Geom. Topol. \textbf{3} (2003), 569--586.

\bibitem{gay}
D.~Gay, \emph{Four-dimensional symplectic cobordisms containing three-handles}, 
Geom. Topol. \textbf{10} (2006) 1749--1759

\bibitem{giroux}
E. Giroux, \emph{Structures de contact en dimension trois et bifurcations des feuilletages de surfaces},
Invent. Math. \textbf{141} (2000), 615--689.

\bibitem{gi}
E. Giroux, \emph{G\'{e}ometrie de contact: de la dimension trois
vers les dimensions sup\'{e}rieures,} Proceedings of the
International Congress of Mathematicians (Beijing 2002), Vol. II,
405--414.

\bibitem{gs}
R.~Gompf and A.~Stipsicz, {\em 4--manifolds and Kirby calculus},
Grad. Stud. Math., Vol. {\bf 20}, AMS, 1999.

\bibitem{h}
K. Honda, \emph{On the classification of tight contact structures
II}, J. Differential Geom. \textbf{55} (2000), 83--143.

\bibitem{hkm}
K. Honda, W. Kazez and G. Matic, \emph{Right-veering diffeomorphisms
of compact surfaces with boundary I}, preprint,
arXiv:math.GT/0510639.

\bibitem{km}
R.~Kirby and P.~Melvin, \emph{Dedekind sums, $\mu$-invariants and
the signature cocyle}, Math. Ann. \textbf{299} (1994), 231--267.

\bibitem{ko}
M.~Korkmaz and B.~Ozbagci, \emph{On sections of elliptic fibrations
}, Michigan Math. J. (to appear), arXiv:math.GT/0604516.

\bibitem{marti}
J. Martinet, \emph{Formes de contact sur les vari\'et\'es de
dimension 3,} in: Proc. Liverpool Singularity Sympos. II, Lecture
Notes in Math. \textbf{209,} Springer, Berlin (1971), 142--163.



\bibitem{oss}
P. Ozsv\'{a}th, A. Stipsicz and Z. Szab\'{o}, {\em Planar open books
and Floer homology},  Int. Math. Res. Not. (2005), 3385--3401.

\bibitem{rol}
D. Rolfsen, \emph{Knots and links,} Publish or Perish, 1976.



\bibitem{tw}
W.P. Thurston and H.E. Winkelnkemper, \emph{On the existence of
contact forms}, Proc. Amer. Math. Soc. \textbf{52} (1975), 345--347.

\bibitem{vh}
J.~Van Horn, \emph{} Ph.D. Dissertation, University of Texas at
Austin, in preparation.


\end{thebibliography}

\providecommand{\bysame}{\leavevmode\hbox
to3em{\hrulefill}\thinspace}

\end{document}